\documentclass[12pt,a4paper]{article}
\usepackage{amsmath,amsthm,amssymb,txfonts,bm}
\usepackage[utf8]{inputenc}

\theoremstyle{definition}
\newtheorem{definition}{Definition}
\newtheorem{thm}[definition]{Theorem}
\newtheorem{lem}[definition]{Lemma}

\newtheorem{remark}[definition]{Remark}

\def\F#1#2#3#4{{}_{2}F_{1}\biggl(\genfrac..{0pt}{}{{#1},\,{#2}}{#3};#4\biggr)}
\def\p#1#2#3#4{{}_{2}\phi_{1}\biggl(\genfrac..{0pt}{}{{#1},\,{#2}}{#3};q,\,#4\biggr)}
\def\pp#1#2#3#4#5#6#7#8{{}_{4}\phi_{3}\biggl(\genfrac..{0pt}{}{{#1},\,{#2},\,{#3},\,{#4}}{{#5},\,{#6},\,{#7}};q,\,#8\biggr)}
\def\Q#1#2#3#4#5#6#7#8{Q\biggl(\genfrac..{0pt}{}{{#1},\,{#2}}{#3};{#4}\,;\genfrac..{0pt}{}{{#5},\,{#6}}{#7};#8\biggr)}
\def\R#1#2#3#4#5#6#7#8{R\biggl(\genfrac..{0pt}{}{{#1},\,{#2}}{#3};{#4}\,;\genfrac..{0pt}{}{{#5},\,{#6}}{#7};#8\biggr)}

\makeatletter
\@addtoreset{equation}{section}

\makeatother

\begin{document}

\title{A $q$-analogue of Gosper's strange evaluation \\of the hypergeometric series}
\author{Yuka Yamaguchi}
\date{\today}

\maketitle

\begin{abstract}
In 1977, Gosper conjectured many strange evaluations of hypergeometric series. 
One of them is a ${}_{2}F_{1}$-series identity with two free parameters, which was proved by Ebisu (2013), Chu (2017), and Campbell (2023) in different ways.
In this paper, we present a $q$-analogue of the ${}_{2}F_{1}$-series identity,  along with its generalization, by using three-term relations for the ${}_{2}\phi_{1}$ basic hypergeometric series. 


\end{abstract}

\section{Introduction}
In 1977, Gosper conjectured many strange evaluations of hypergeometric series. 
One of them is the following identity for the ${}_{2}F_{1}$ hypergeometric series with two free parameters: 
\begin{align}\label{Gosper}
\F{1 - \alpha}{\beta}{\beta + 2}{\frac{\beta}{\alpha + \beta}} = (\beta + 1) \left(\frac{\alpha}{\alpha + \beta} \right)^\alpha, 
\end{align}
where 
\begin{align*}
\F{\alpha}{\beta}{\gamma}{x} 
= {}_{2}F_{1} (\alpha, \beta; \gamma; x) 
:=  \sum_{n = 0}^{\infty} 
\frac{(\alpha)_{n} (\beta)_{n}}{(1)_{n} (\gamma)_{n}} x^n
\end{align*}
with $(\alpha)_{n} := \Gamma(\alpha + n) / \Gamma(\alpha)$. 
The identity~($\ref{Gosper}$) was proved by Ebisu \cite{Eb1}, Chu \cite{Chu}, and Campbell \cite{Campbell} in different ways. 
In this paper, we give a $q$-analogue of ($\ref{Gosper}$), along with its generalization. 
Ebisu's proof~\cite{Eb1} of ($\ref{Gosper}$) is based on the results of \cite{Eb0} regarding three-term relations for the ${}_{2}F_{1}$-series. 
Our approach to providing a $q$-analogue of ($\ref{Gosper}$) is based on the results of \cite{Y1} regarding three-term relations for the ${}_{2}\phi_{1}$ basic hypergeometric series, which is a $q$-analogue of the ${}_{2}F_{1}$-series. 
Although our proof is partially similar to Ebisu's proof, it is not just an analogy: Ebisu employed an expression for ${}_{2}F_{1}(\alpha, 1; \gamma; x)$ in terms of the incomplete beta function, whereas we only use three-term relations. 

The ${}_{r+1} \phi_{r}$ basic hypergeometric series is defined by 
\begin{align*}
{}_{r+1}\phi_{r}\biggl(\genfrac..{0pt}{}{a,\, b_{1}, \dotsc, b_{r}}{c_{1}, \dotsc, c_{r}}; q,\, x\biggr) 
&= {}_{r+1}\phi_{r} (a, b_{1}, \dotsc, b_{r}; c_{1}, \dotsc, c_{r}; q, x) \\
&:= \sum_{n = 0}^{\infty} 
\frac{(a;q)_{n} (b_{1};q)_{n} \dotsm (b_{r};q)_{n}}{(q;q)_{n} (c_{1};q)_{n} \dotsm (c_{r};q)_{n}} x^{n}, \nonumber 
\end{align*}
where $(a; q)_{n} := (a; q)_{\infty} / (a q^{n}; q)_{\infty}$ with 
$(a; q)_{\infty} := \prod_{j = 0}^{\infty} (1 - a q^{j})$ and it is assumed that $\lvert q \rvert < 1$. 
In the limit as $q \to 1$, the series ${}_{2}\phi_{1} (q^\alpha, q^\beta; q^\gamma; x)$ becomes ${}_{2}F_{1} (\alpha, \beta; \gamma; x)$. 

Our main result is the following. 
\begin{thm}\label{thm:1}
Let $m$ be an integer with $m \geq 2$. 
We assume that $c,\, a/c \notin q^{\mathbb{Z}} \cup \{0 \}$. 
Then, for any root $\lambda$ of the polynomial ${}_{2}\phi_{1} (q/a, q^{1 - m}; q^2/c; q, a q^m x / c)$ in the variable $x$ of degree at most $m - 1$, the following hold: 
\begin{align}
\p{a}{q^m}{c}{\lambda} 
&= \lambda^{1-m} f(\lambda), \label{(1.2)}\\
\p{c/a}{c q^{-m}}{c}{\frac{a q^m \lambda}{c}} 
&= \frac{(\lambda; q)_{\infty}}{(a q^m \lambda/ c; q)_{\infty}} \lambda^{1-m} f(\lambda), \label{(1.3)}
\end{align}
where $f(x)$ is the polynomial in $x$ of degree at most $m - 2$ defined by 
\begin{align*}
f(x) 
&:= - \frac{(q/c;q)_{m-1} (cq^{-m};q)_{m-1}}{(aq/c;q)_{m-1} (q;q)_{m-1}} q^{m-1} \\
&\quad \times \sum_{j=0}^{m-2} \frac{(q/a;q)_{j} (q^{1-m};q)_{j}}{(q;q)_{j} (q^2/c;q)_{j}} \left(\frac{aq^m x}{c} \right)^{j} 
\pp{q^{-j}}{q}{a}{cq^{-j-1}}{q^{m-j}}{aq^{-j}}{cq^{1-m}}{q}. 
\end{align*}
When $m = 2$, we obtain 
\begin{align}
\p{a}{q^2}{c}{\frac{q^2-c}{(q-a)q}} 
&= \frac{(q-a) (q-c)}{(1-q) (c-aq)}, \label{(1.4)} \\
\p{c/a}{c q^{-2}}{c}{\frac{(q^2-c)aq}{(q-a)c}} 
&= \frac{(\lambda; q)_{\infty}}{(a q^2 \lambda/ c; q)_{\infty}} 
\frac{(q-a) (q-c)}{(1-q) (c-aq)}, \label{(1.5)} 
\end{align}
where $\lambda = (q^2-c) / \{(q-a)q \}$. 
\end{thm}

\begin{remark}
The degree of ${}_{2}\phi_{1} (q/a, q^{1 - m}; q^2/c; q, a q^m x / c)$ is exactly $m-1$ if we further assume that $a \neq q^{i}$ ($1 \leq i \leq m-1$). 
\end{remark}

From the $q$-binomial and the binomial theorems, it holds that  
\begin{align*}
\frac{(q^{\alpha}z;q)_{\infty}}{(z;q)_{\infty}} 
= \sum_{n = 0}^{\infty} \frac{(q^{\alpha};q)_{n}}{(q;q)_{n}} z^n 
\; \xrightarrow{q \to 1} \;
\sum_{n = 0}^{\infty} \frac{(\alpha)_{n}}{(1)_{n}} z^n 
= (1-z)^{-\alpha}, 
\end{align*}
where $\lvert z \rvert < 1$ (see \cite[(1.3.2), (1.3.1)]{GR}). 
Thus, setting $(a, c) = (q^{\alpha+\beta+1}, q^{\beta+2})$ in $(\ref{(1.5)})$ yields a $q$-analogue of $(\ref{Gosper})$. 
Also, $(\ref{(1.4)})$ and $(\ref{(1.5)})$ with $(a, c) = (q^\alpha, q^\gamma)$ are $q$-analogues of \cite[(1.5), (1.6)]{Eb1}, respectively. 
Moreover, $q$-analogues of \cite[(1.2), (1.3)]{Eb1}, which are alternative expressions of $(\ref{(1.2)})$ and $(\ref{(1.3)})$, are given in Section~$3$ (see Remark~$\ref{remark}$). 


\section{Three-term relations}
In this section, we introduce three-term relations for the ${}_{2}\phi_{1}$-series and give two lemmas that are used to prove Theorem~$\ref{thm:1}$. 

Any three ${}_{2}\phi_{1}$-series, where respective parameters $a, b, c$ and a variable $x$ of them are shifted by integer powers of $q$, satisfy a linear relation with coefficients that are rational functions of $a, b, c, q$, and $x$. We call such a relation a three-term relation for the ${}_{2}\phi_{1}$-series.  
In \cite[Theorem~$1$]{Y1}, for any integers $k, l, m$, and $n$, 
the coefficients $Q$ and $R$ of the three-term relation 
\begin{align}\label{3tr}
\p{a q^{k}}{b q^{l}}{c q^{m}}{x q^{n}} 
= Q \cdot \p{a q}{b q}{c q}{x} + R \cdot \p{a}{b}{c}{x}
\end{align}
are explicitly expressed. 
Also, expressions for $Q$ and $R$ as sums of products of two ${}_{2}\phi_{1}$-series are provided in \cite[Lemma~$15$]{Y1}. 
(The $n = 0$ cases are also available in \cite{Y0}). 
In accordance with the assumption made in \cite{Y1}, we also assume that 
\begin{align}\label{assumption}
a, b, c, \frac{a}{b}, \frac{c}{a}, \frac{c}{b} \notin q^{\mathbb{Z}} \cup\{0 \}. 
\end{align}
We denote the coefficients of $(\ref{3tr})$ by 
\begin{align*}
Q = \Q{k}{l}{m}{n}{a}{b}{c}{x}, \quad 
R = \R{k}{l}{m}{n}{a}{b}{c}{x}. 
\end{align*}
From \cite[Theorem~$1$, Lemma~$15$]{Y1}, the following holds: 
When $k \leq l$ and $0 \leq m \leq k+l+n$, the coefficient $Q$ of $(\ref{3tr})$ can expressed by 
\begin{align}\label{(2.3)}
&\Q{k}{l}{m}{n}{a}{b}{c}{x} \nonumber \\
&= - \frac{(1-a) c}{(q-c) (1-c)} \frac{(x;q)_{n}\, x^{1-m}}{(abqx/c;q)_{k+l-m+n-1}} \nonumber \\
&\;\times \left\{\frac{(aq/c;q)_{k-m} (bq/c;q)_{l-m} (c;q)_{m}}{(q^2/c;q)_{-m} (a;q)_{k} (bq;q)_{l-1}} (cq^{m-1})^{-n} \right. \nonumber \\
&\qquad \times \p{q^{1-k}/a}{q^{1-l}/b}{q^{2-m}/c}{\frac{abq^{k+l-m+n}}{c}x} 
\p{a}{b}{c}{x} \nonumber \\
&\qquad \left. - (1-b) x^{m} 
\p{cq^{m-k}/a}{cq^{m-l}/b}{cq^m}{\frac{abq^{k+l-m+n}}{c}x} 
\p{aq/c}{bq/c}{q^2/c}{x} \right\}. 
\end{align}
Also, from \cite[Theorem~$1$]{Y1}, the following holds: When $(k, l, m, n) = (1, m, m, 0)$ with $m \geq 2$, the coefficient $R$ of $(\ref{3tr})$ can expressed by 
\begin{align}\label{(2.4)}
&\R{1}{m}{m}{0}{a}{b}{c}{x} \nonumber \\
&= \frac{(q^{1-m}/c;q)_{m} (cq;q)_{m-2}}{(aq^{2-m}/c;q)_{m-1} (bq;q)_{m-1}} 
c q^{m-1} x^{1-m} \nonumber \\
&\quad \times \sum_{j=0}^{m-2} \frac{(1/a;q)_{j} (q^{1-m}/b;q)_{j}}{(q;q)_j (q^{2-m}/c;q)_j} 
\left(\frac{abqx}{c} \right)^j 
\pp{q^{-j}}{cq^{m-j-1}}{aq}{bq}{cq}{aq^{1-j}}{bq^{m-j}}{q}, 
\end{align}
where we used the fact that $B_{j}$ in \cite[Theorem~$1$]{Y1} equals $0$ for $j < 0$. 
As for $(\ref{(2.3)})$ and $(\ref{(2.4)})$, by analytic continuation, we can drop the assumption that $b \neq q^{0}$ made in $(\ref{assumption})$, 
and therefore we obtain the following two lemmas. 

\begin{lem}\label{lem:3}
When $k \leq l$ and $0 \leq m \leq k+l+n$, we have 
\begin{align*}
&\Q{k}{l}{m}{n}{a}{1}{c}{x} \\
&= - \frac{(1-a) c}{(q-c) (1-c)} \frac{(x;q)_{n}\, x^{1-m}}{(aqx/c;q)_{k+l-m+n-1}} \\
&\quad \times \frac{(aq/c;q)_{k-m} (q/c;q)_{l-m} (c;q)_{m}}{(q^2/c;q)_{-m} (a;q)_{k} (q;q)_{l-1}} (cq^{m-1})^{-n}
\p{q^{1-k}/a}{q^{1-l}}{q^{2-m}/c}{\frac{aq^{k+l-m+n}}{c}x}. 
\end{align*}
\end{lem}

\begin{lem}\label{lem:4}
When $m \geq 2$, we have
\begin{align*}
&\R{1}{m}{m}{0}{a}{1}{c}{x} \\
&= \frac{(q^{1-m}/c;q)_{m} (cq;q)_{m-2}}{(aq^{2-m}/c;q)_{m-1} (q;q)_{m-1}} 
c q^{m-1} x^{1-m} \nonumber \\
&\quad \times \sum_{j=0}^{m-2} \frac{(1/a;q)_{j} (q^{1-m};q)_{j}}{(q;q)_j (q^{2-m}/c;q)_j} 
\left(\frac{aqx}{c} \right)^j 
\pp{q^{-j}}{cq^{m-j-1}}{aq}{q}{cq}{aq^{1-j}}{q^{m-j}}{q}. 
\end{align*}
\end{lem}

\section{Proof of Theorem~$\mathbf{\ref{thm:1}}$}
We prove Theorem~$\ref{thm:1}$. 
When $(k, l, m, n) = (1, m, m, 0)$ with $m \geq 2$, by substituting $b = 1$ into $(\ref{3tr})$, we have
\begin{align}\label{(3.1)}
\p{a q}{q^{m}}{c q^{m}}{x} 
= \Q{1}{m}{m}{0}{a}{1}{c}{x} \cdot \p{a q}{q}{c q}{x} 
+ \R{1}{m}{m}{0}{a}{1}{c}{x}.  
\end{align}
Let $x = \lambda$ be a root of the polynomial ${}_{2}\phi_{1}(1/a, q^{1-m}; q^{2-m}/c; q, aqx/c)$ in $x$. 
Then, by Lemmas~$\ref{lem:3}$ and $\ref{lem:4}$, substituting $x = \lambda$ into $(\ref{(3.1)})$ yields 
\begin{align*}
\p{a q}{q^{m}}{c q^{m}}{\lambda} 
&= \R{1}{m}{m}{0}{a}{1}{c}{\lambda} \\
&= \frac{(q^{1-m}/c;q)_{m} (cq;q)_{m-2}}{(aq^{2-m}/c;q)_{m-1} (q;q)_{m-1}} 
c q^{m-1} \lambda^{1-m} \\
&\;\times \sum_{j=0}^{m-2} \frac{(1/a;q)_{j} (q^{1-m};q)_{j}}{(q;q)_j (q^{2-m}/c;q)_j} 
\left(\frac{aq\lambda}{c} \right)^j 
\pp{q^{-j}}{cq^{m-j-1}}{aq}{q}{cq}{aq^{1-j}}{q^{m-j}}{q}. 
\end{align*}
By replacing $(a, c)$ with $(aq^{-1}, cq^{-m})$, and by rewriting $(q/c;q)_m (cq^{1-m};q)_{m-2}\,cq^{-1}$ as $-(q/c;q)_{m-1} (cq^{-m};q)_{m-1}\,q^{m-1}$, we obtain $(\ref{(1.2)})$. 
Also, applying Heine's \cite{Heine} (also available in \cite[(1.4.3)]{GR}) transformation formula 
\begin{align}\label{Heine}
\p{a}{b}{c}{x} = \frac{(abx/c;q)_\infty}{(x;q)_\infty} \p{c/a}{c/b}{c}{\frac{abx}{c}}
\end{align}
to the left-hand side of $(\ref{(1.2)})$ yields $(\ref{(1.3)})$. 

\begin{remark}\label{remark}
In the above, we set $(k, l, m, n) = (1, m, m, 0)$ with $m \geq 2$ to derive $(\ref{(1.2)})$. 
Alternative expressions of $(\ref{(1.2)})$ can be obtained by choosing other values of $(k, l, m, n)$ such that $l \geq \max \{k, 2 \}$, $0 \leq m \leq k+l$, and $n=0$. 
We now fix such a choice of $(k, l, m, n)$. 
Then, by Lemma~$\ref{lem:3}$, 
substituting $b=1$ and $x=\lambda$ into $(\ref{3tr})$, where $\lambda$ is a root of the polynomial ${}_{2}\phi_{1}(q^{1-k}/a, q^{1-l}; q^{2-m}/c; q, aq^{k+l-m}x/c)$ in $x$, yield  
\begin{align*}
\p{a q^k}{q^l}{c q^m}{\lambda} 
= \R{k}{l}{m}{0}{a}{1}{c}{\lambda}. 
\end{align*}
Replacing $(a, c)$ with $(aq^{-k}, cq^{-m})$ gives
\begin{align*}
\p{a}{q^{l}}{c}{\lambda} 
= \R{k}{l}{m}{0}{a q^{-k}}{1}{c q^{-m}}{\lambda}, 
\end{align*}
where $\lambda$ is a root of the polynomial ${}_{2}\phi_{1}(q/a, q^{1-l}; q^2/c; q, aq^{l}x/c)$ in $x$. 
Therefore, we can obtain an alternative expression of $(\ref{(1.2)})$ from \cite[Theorem~$1$]{Y1}. 
If we further decide to use \cite[Lemma~$15$]{Y1}, we can express it as a sum of products of two ${}_{2}\phi_{1}$-series. 
In particular, when $(k, l, m, n) = (1, l, 1, 0)$ with $l \geq 2$, 
we obtain 
\begin{align*}
\p{a}{q^{l}}{c}{\lambda} 
&= \frac{1}{(aq\lambda/c;q)_{l-1}} \p{c/a}{cq^{-l}}{c}{\frac{aq^{l}\lambda}{c}} {}_{1}\phi_{0}\biggl(\genfrac..{0pt}{}{aq/c}{\mathrm{-}};q,\,\lambda\biggr) \nonumber \\
&\quad - \frac{(q^2/c;q)_{l-1}}{(aq\lambda/c;q)_{l-1} (q;q)_{l-1}} 
\p{q/a}{q^{1-l}}{q^2/c}{\frac{aq^{l}\lambda}{c}} \p{a}{q}{c}{\lambda}. 
\end{align*}
This is a $q$-analogue of \cite[(1.2)]{Eb1}. 
Also, applying $(\ref{Heine})$ to the left-hand side yields a $q$-analogue of \cite[(1.3)]{Eb1}. 
\end{remark}

\clearpage

\section*{Acknowledgment}
This work was supported by JSPS KAKENHI Grant Number JP25KJ0266. 

\bibliography{reference}

\medskip
\begin{flushleft}
Faculty of Education \\ 
University of Miyazaki \\ 
1-1 Gakuen Kibanadai-nishi \\ 
Miyazaki 889-2192 Japan \\ 
{\it Email address}: y-yamaguchi@cc.miyazaki-u.ac.jp 
\end{flushleft}

\end{document}